\documentclass[12pt]{amsart}

\usepackage{amsmath,amsfonts,amsthm,amsopn,cite,mathrsfs}
\usepackage{epsfig,verbatim}
\usepackage{subfigure}
\newcommand{\tfa}{time-frequency analysis}

\newcommand{\stft}{short-time Fourier transform}

\newcommand{\tf}{time-frequency}

\newcommand{\fif}{if and only if}
\newcommand{\tfs}{time-frequency shift}

\newcommand{\modsp}{modulation space}

\newtheorem{tm}{Theorem}[section]
\newtheorem{lemma}[tm]{Lemma}
\newtheorem{prop}[tm]{Proposition}
\newtheorem{cor}[tm]{Corollary}

\newcommand{\rem}{\noindent\textsl{REMARK:}}

 \theoremstyle{definition}
 \newtheorem{definition}{Definition}

\newcommand{\beqa}{\begin{eqnarray*}}
\newcommand{\eeqa}{\end{eqnarray*}}

\newcommand{\field}[1]{\mathbb{#1}}
\newcommand{\bR}{\field{R}}        
\newcommand{\bN}{\field{N}}        
\newcommand{\bZ}{\field{Z}}        
\newcommand{\bC}{\field{C}}        
\newcommand{\bT}{\field{T}}        %
 \def\cS{\mathcal{S}}

 \def\cB{\mathcal{B}}
 
 \def\cG{\mathcal{G}}

\def\rd{\bR^d}

\def\rdd{{\bR^{2d}}}
\def\zdd{{\bZ^{2d}}}

\def\lrd{L^2(\rd)}

\def\zd{\bZ^d}

\def\intrd{\int_{\rd}}
\def\intrdd{\int_{\rdd}}

\def\<{\left<}
\def\>{\right>}

\def\inv{^{-1}}

\def\mv1{M_v^1}

\def\phas{(x,\xi )}

\renewcommand{\qed}{\hfill \rule{7pt}{8pt} \vskip .2truein}   

\newcommand{\lone}{\ell ^1}

\newcommand{\duall}{\Lambda ^{\circ}}
\newcommand{\latt}{\Lambda}
\newcommand{\gab}{\cG (g, \latt )}
\newcommand{\gabdual}{\cG (g, \duall )}
\newcommand{\vf}{\varphi}
\newcommand{\bba}{\textbf{a}}
\newcommand{\bbb}{\textbf{b}}
\newcommand{\bbc}{\textbf{c}}
\newcommand{\imp}{$\, \Rightarrow \, $}
\newcommand{\twi}{\, \natural \, }
\begin{document}
\begin{abstract}
We prove fourteen equivalent conditions for a set of \tfs s on a
lattice $\Lambda $, $\{
e^{2\pi i \lambda _2\cdot t} g(t-\lambda _1) : (\lambda _1, \lambda _2) \in
\Lambda \}\subseteq \lrd$, to be a frame for $\lrd $. Remarkably,
several of these conditions 
can  be formulated without an inequality. In particular, instead of
checking the invertibility of the frame operator on $\lrd $,
it suffices to verify that it is one-to-one on a certain subspace of
tempered distributions.  
\end{abstract}

\title{Gabor Frames without Inequalities}
\author{Karlheinz Gr\"ochenig}
\address{Faculty of Mathematics \\
University of Vienna \\
Nordbergstrasse 15 \\
A-1090 Vienna, Austria}
\email{karlheinz.groechenig@univie.ac.at}
\subjclass{}
\date{March 13, 2007}
\keywords{Gabor frame, \modsp , twisted convolution, injectivity
  implies invertibility}
\thanks{K.~G.~was supported by the Marie-Curie Excellence Grant
  MEXT-CT 2004-517154}
\maketitle

\section{Introduction}

   An operator on a finite-dimensional vector space is
invertible, \fif\ it is one-to-one.  On an  infinite-dimensional vector
space, this characterization is hardly ever true and amounts to a
mathematical miracle that says something deep about the operator. 

In this paper we investigate the principle ``injectivity implies
invertibility'' in the context of \tfa\ (phase-space analysis).  We
show that several operators associated to a Gabor frame are invertible
on Hilbert space, \fif\ they are one-to-one on a larger space. This
insight seems to be completely new and is 
rather surprising. 

Gabor frames are an important and convenient tool for phase-space
expansions and the  \tfa\ of distributions. The theory of Gabor frames
is well developed, as is seen from the
textbooks~\cite{chr03,daubechies92,book}. Several explicit
constructions  and many 
general characterizations are known. 
In all these characterizations the difficulty is to 
check the \emph{invertibility} of some  operator on a Hilbert space: this  is 
either  the associated  frame operator on $\lrd $, or  some Gramian
matrix on $\ell ^2(\zdd )$, or even
of a whole family of matrices acting on $\ell ^2(\zd )$ . 
 Equivalently, the invertibility amounts to finding a positive
 \emph{lower} bound for  an \emph{inequality}, which is inevitably a
 hard problem.  

Gabor frames are also implicit in Rieffel's work on projective modules
over non-commutative tori~\cite{rieffel88} and thus play a (not yet
fully understood) role  in non-commutative geometry (see Luef's work on
the connection between the two fields~\cite{lue06}).

 We offer a number  of new characterizations of Gabor
frames that do not require inequalities or invertibility.   In each
case  it suffices to verify the injectivity of some 
operators associated to a Gabor frame on a larger space instead
proving its invertibility on a Hilbert space. These results
reveal a remarkable  phenomenon in
phase-space analysis and, to us, are completely unexpected.  
The results do not come easily: we will use some of the deepest
results about Gabor frames, such as the duality theorem of Janssen and
Ron-Shen, Wiener's Lemma for twisted convolution and the rotation
algebra.

The paper is organized as follows. In Section~2 we collect the basic 
definitions required to formulate our results on Gabor frames. The
main result is stated and commented in Section~3. It proof is carried
out in Section~4 after a brief summary of several results from
abstract harmonic analysis.  Section~5 provides a few examples and
further perspectives, such as the index of a Gabor system. 

\section{Gabor Frames and Modulation Spaces --- Basic Definitions}

To make our statements precise, let us introduce the main concepts
required to treat Gabor frames. For the detailed exposition of Gabor
frames and \tfa\ we refer to the
books~\cite{chr03,daubechies92,book}. 

Given a point $z= \phas , x,\xi \in \rd, $ in phase space $\rdd $, we
consider the phase-space shifts (\tfs s)  acting on a function $f$
\begin{equation}
  \label{eq:1}
  \pi (z) f(t) = M_\xi T_x f(t) = e^{2\pi i\xi \cdot t} f(t-x), \qquad
  t, x, \xi \in \rd \, .
\end{equation}

A \tf\ lattice $\Lambda $ is a discrete  subgroup of $\rdd $,
of the form $\Lambda = A \zdd $ for some invertible real-valued $2d
\times 2d$-matrix. The \emph{adjoint lattice } is 
$\duall = \{ \mu \in \rdd : \pi (\lambda )
\pi (\mu )  = \pi (\mu ) \pi (\lambda ) \, \text{ for all } \, \lambda
\in \Lambda \}$. If $\Lambda = \alpha \zd \times \beta \zd$, then
$\duall = \beta \inv \zd \times \alpha \inv \zd $. The adjoint lattice
should be distinguished from the 
dual lattice $\Lambda ^\perp = \alpha \inv \zd \times \beta \inv \zd$
that is more common is harmonic analysis. 

Fix a $g\in \lrd $, then the set of   \tfs s with respect to the
lattice $\Lambda $ is denoted by 
\begin{equation}
  \label{eq:2}
  \gab = \{ \pi (\lambda ) g : \lambda \in \Lambda \} \, ,
\end{equation}
and similarly $\gabdual $ is defined as $  \gabdual = \{ \pi (\mu ) g
: \mu  \in \duall \}  $. 
Sets of the form $\gab $ and $\gabdual $ are called \emph{Gabor
  systems}.

Associated to every Gabor system there are four  canonical operators:
the analysis operator,  the synthesis operator, the frame operator, and
the Gramian. Precisely, these operators are defined as follows.

\begin{definition}
Fix a non-zero  $g\in \lrd $ and let  $\Lambda \subseteq \rdd $ be a
lattice. 
The \emph{coefficient operator} $C_{g,\latt }$ associated to a test function
$g$ and a lattice $\latt $ maps functions/distributions
to sequences on $\latt $ and  is defined to be
\begin{equation}
  \label{eq:13}
  (C_{g,\latt } f) (\lambda ) = \langle f, \pi (\lambda ) g \rangle
  \qquad f\in \lrd , \lambda \in \latt \, .
\end{equation}
The \emph{synthesis operator} $D_{g,\latt } $   maps sequences to
functions/distributions and is 
\begin{equation}
  \label{eq:14}
  D_{g,\latt } \mathbf{c} = \sum _{\lambda \in \latt } c_\lambda \pi
  (\lambda ) g \, ,
\end{equation}
whenever the series is defined (for instance,  for finite sequences
$\mathbf{c}$). 

The composition $S_{g,\latt } = D_{g,\latt }  C_{g,\latt } $ is the
\emph{frame operator} corresponding to the Gabor system $\gab $ and
maps functions to functions, when well-defined. 

Finally, the operator $G_{g, \latt }  = C_{g, \latt } D_{g, \latt } $
is the \emph{Gramian operator} mapping sequences indexed by $\latt $ to
sequences. Viewed as a matrix, $G_{g, \latt } $ has the entries
$G_{\lambda , \lambda '} = \langle \pi (\lambda ' ) g, \pi (\lambda )
g\rangle $ for $\lambda , \lambda ' \in \latt $. 
\end{definition}

  \begin{definition}
    The set $\gab $ is called a \emph{Gabor frame} (or Weyl-Heisenberg
    frame), if  there exist constants $A,B >0$, such that 
      \begin{equation}
        \label{eq:3}
        \qquad A \|f\|_2^2 \leq  \sum _{\lambda \in \Lambda } |\langle f, \pi
        (\lambda ) g \rangle |^2 \leq B\|f \|_2^2 \qquad \text{for all
        } \, f \in \lrd \, .
      \end{equation}
Further,     the set $\gab $ is called a \emph{Riesz sequence}, if there
    exist constants $A', B' > 0$ such that 
    \begin{equation}
      \label{eq:4}
      A' \|\mathbf{c}\|_2 \leq \|\sum _{\lambda \in \latt } c_\lambda \pi
      (\lambda  ) g \|_2 \leq B' \|\mathbf{c}\|_2 
    \end{equation}
holds for all finite sequence $\mathbf{c}$. 
  \end{definition}
Since $\langle Sf,f\rangle =  \sum _{\lambda \in \Lambda } |\langle f, \pi
        (\lambda ) g \rangle |^2$, the inequalities in  \eqref{eq:3}
        express  that the Gabor  frame operator $S_{g,\latt} $ is
      bounded and   invertible on $\lrd $.  To prove that $\cG (g, \latt )$ is a
        frame, it is therefore necessary to either prove the
        invertibility of the operator $S_{g, \latt }$  or to prove the
        inequalities \eqref{eq:3}. Whereas the upper inequality is
        easy (it amounts to the boundedness of $S_{g,\latt }$), the
        lower  lower inequality is difficult (it amounts to the
        invertibility of $S_{g,\latt }$). Several criteria for Gabor frames have been
found, notably the Wexler-Raz conditions~\cite{DLL95,janssen95}, the
Ron-Shen duality~\cite{ron-shen97}, the characterization by the Ron-Shen
matrices~\cite{ron-shen97,book}.  

Likewise, $\cG (g, \Lambda )$ is a Riesz sequence, \fif\ the  Gramian
is invertible on $\ell ^2 (\zd )$.  This follows from the identity  
\begin{eqnarray*}
\|\sum _{\lambda \in \latt } c_\lambda \pi       (\lambda ) g \|_2^2 &=& \sum
_{\lambda, \lambda '  \in \latt } c_\lambda \overline{c_{\lambda '}} \langle  \pi
      (\lambda ) g, \pi (\lambda ')g\rangle \\
&=&   \sum _{\lambda, \lambda '  \in \latt } c_\lambda \overline{c_{\lambda '}}
G_{\lambda ' ,\lambda } = \langle G_{g,\latt }  \mathbf{c}, \mathbf{c} \rangle  \, .
\end{eqnarray*}

We note that whenever a duality $\langle \cdot , \cdot \rangle $ is
suitably defined and extends the inner product on $\lrd $ (in
particular, it is conjugate linear in the second term!), then the
synthesis operator is adjoint to the analysis operator and vice
versa, informally 
\begin{equation}
  \label{eq:8}
  C_{g,\latt } ^* = D_{g, \latt } \, .
\end{equation}

\textbf{Modulation Spaces.}  

\begin{definition}
Fix a non-zero Schwartz function $\vf $
(preferrably the Gaussian $\vf (t) = e^{-\pi t\cdot t}$). We say that a
distribution $f\in \cS ' (\rd ) $ belongs to the \emph{\modsp }~  $M^p
= M^p(\rd )$, if 
\begin{equation}
  \label{eq:6}
  \|f\|_{M^p} := \Big(\intrdd | \langle f, \pi (z) \vf \rangle |^p \,
  dz \Big)^{1/p} = \|V_\vf f\|_p  <
  \infty \, .
\end{equation}
\end{definition}

 Each $M^p(\rd ), 1\leq p \leq \infty, $ is a  Banach space, and
the  definition is 
independent of the test function $\vf \in \cS (\rd
)$.
If $1\leq p< \infty $, then the dual space of $M^p$ is $M^{p'}$ where
$1/p+1/p' = 1$~\cite[Thm.~11.3.6]{book}. We note that for $p=2$ we
obtain $M^2 = \lrd $.  
See~\cite{feichtinger81,feichtinger89,feichtinger80cras,feiSTSIP,book}
for the  many beautiful  
properties of the \modsp s and their many generalizations. 

The transform 
$$
\langle f, \pi (z) \vf \rangle = \langle f, M_\xi T_x \vf \rangle =
\intrd f(t) \, \overline{g(t-x)} \, e^{-2\pi i x\cdot t} \, dt
$$
is the so-called \stft\ (also called Gabor transform, ambiguity
function, coherent state transform). It measures the phase-space (\tf )
content at $z= (x,\xi )$ in phase space. So $M^1$ consists of all
those $L^2$-functions whose \stft\ is absolutely integrable, and
$M^\infty $ contains exactly those tempered distributions with bounded
\stft . It is not hard to see that the Schwartz class is a subspace of
$ M^1$. In \tfa\ and phase-space
analysis, $M^1$ is therefore often used as  an appropriate space of
test functions, 
and $M^\infty $ serves as a suitable space of
distributions~\cite{feichtinger80cras,FZ98}.     The \modsp s are  
tailored to the needs of \tfa , and they arise inevitably, whenever
a problem involves the \tfs s $\pi (z)$. For an exposition from
scratch see~\cite[Ch.~11-13]{book}, for a detailed history and a
  comprehensive annotated list of references see Feichtinger's
  beautiful article~\cite{feiSTSIP}.



 With all definitions in place, we  verify first when the operators
 associated to a Gabor system are bounded. 

\begin{lemma}\label{bound}
Assume that $g\in M^1, g\neq 0$. Then
\begin{itemize}
\item[(i)] $C_{g,\latt }$ maps $M^p(\rd )$ into $\ell ^p (\Lambda )$
  and $
\|C_{g,\latt } f \|_p \leq C \|g\|_{M^1} \|f\|_{M^p} $. The constant
$C$ depends only on the lattice, but not on  $g$ and $f$. 
\item[(ii)] $D_{g,\latt }$ maps $\ell ^p (\Lambda )$ into $M^p(\rd )$
  and $\|D_{g,\latt } \mathbf{c} \|_p \leq C \|g\|_{M^1} \|\mathbf{c}\|_{\ell
^p} $.
\item[(iii)] The frame operator $S_{g,\latt }$ maps $M^p(\rd )$ into
  $M ^p (\rd )$  and \\  $\|S_{g,\latt } f \|_p \leq C^2 \|g\|_{M^1}^2
  \|f\|_{M^p} $. 
\item[(iv)] The Gramian $G_{g,\latt }$ maps $\ell ^p(\latt  )$ into
  $\ell ^p (\Lambda )$   and $\|G_{g,\latt } \bbc  \|_p \leq C^2 \|g\|_{M^1}^2 \|\bbc\|_{p} 
$.
\end{itemize}
  \end{lemma}

These statements are well known and can be found in various sources,
see~\cite{fg97jfa} and~\cite[Cor.~12.1.12]{book} for the proofs and
detailed references. 

For later use, we note that item (i) with $f=g$ implies that for $g\in
M^1$ we have
\begin{equation}
  \label{eq:35}
  \sum _{\lambda \in \latt } |\langle g, \pi (\lambda ) g\rangle | =
  \|C_{g,\latt } g\|_1 \leq C \|g\|_{M_1}^2 \, .
\end{equation}

From (iii) and  (iv) we see that the frame operator $S_{g,\latt} $ and
the Gramian $G_{g,\latt }$ are always bounded. So the right-hand inequalities
in~\eqref{eq:3} and~\eqref{eq:4} are always satisfied when $g\in
M^1$. 


To make sense of infinite series of \tfs s, we need the following lemma. 

\begin{lemma}\label{synthesis}
  If $g\in M^1(\rd )$ and $\mathbf{c} \in \ell ^\infty (\latt )$, then
  the operator $\sum _{\lambda \in \Lambda } c_\lambda \pi (\lambda )$
  is bounded from $M^1(\rd )$ to $M^\infty (\rd )$, and the sum
  converges unconditionally in the weak operator topology. 
\end{lemma}

\begin{proof}
By Lemma~\ref{bound} we have 
  \begin{equation}
    \label{eq:15}
\|C_{g,\latt}f \|_1=     \sum _{\lambda \in \Lambda } |\langle f,
  \pi (\lambda )g \rangle| \leq C\,  \|f\|_{M^1} \,  \|g\|_{M^1} \,
  \qquad \forall f,g\in M^1 \, .
  \end{equation}
Then  we find that 
\begin{eqnarray}
  |\langle \sum _{\lambda \in \Lambda } c_\lambda \pi (\lambda ) f, g
  \rangle | &=& |\sum _{\lambda \in \Lambda } c_\lambda \langle  \pi
  (\lambda ) f, g   \rangle | \notag \\
&\leq & \|\mathbf{c} \|_\infty \, \sum _{\lambda \in \Lambda } |\langle  \pi
  (\lambda ) f, g   \rangle | 
\leq C\, \|\mathbf{c}\|_\infty \, \|f\|_{M^1} \,  \|g\|_{M^1} \, .  \label{eq:36}
\end{eqnarray}
Since  $M^1$ and $M^\infty $ are dual to each other, this inequality
implies  that
 $$\|\sum _{\lambda \in \Lambda } c_\lambda \pi (\lambda
)f\|_{M^\infty } \leq C \| f\|_{M^1} \qquad \text{ for all} \, f\in
M^1 \, .
$$
The weak unconditional convergence of $\sum _\lambda c_\lambda \pi
(\lambda )$   follows immediately from~\eqref{eq:36}. 
\end{proof}

Finally we need a strong form of linear independence of \tfs s. 

\begin{prop}\label{linind}
If $\sum _{\lambda \in \latt } c_\lambda \pi (\lambda ) = 0 $ for some
$\bbc \in \ell ^\infty (\latt )$, then $\bbc = 0$. 

\end{prop}

\begin{proof}
By assumption  we have, for all $g,h\in M^1$ and $z \in \rdd $,
  $$
\sum _{\lambda \in \Lambda } c_\lambda \,  \langle \pi (\lambda ) \pi
(z) g, \pi (z) h \rangle = 0 \, .$$ 
Now $\pi (z) \inv \pi (\lambda ) \pi (z) = e^{2\pi i [ z,\lambda ]}
\pi (\lambda )$, where $[ z,\lambda] =   z_1 \lambda _2 -  z_2
\lambda _1 $ is the symplectic form on $\rdd $. This implies that 
\begin{equation}
  \label{eq:1new}
\sum _{\lambda \in \Lambda } c_\lambda \, \langle \pi (\lambda )  g,  h
\rangle \, \,  e^{2\pi i  \lambda _2 z_1 -  \lambda _1 z_2
  } = 0 \,   
\end{equation}
for all $z\in \rdd $ and all $g,h \in M^1$. 

Equation~\eqref{eq:1new} is an absolutely converging Fourier series on
$\rdd / \Lambda $. Since   it vanishes everywhere, we must have 
$$
 c_\lambda \, \langle \pi (\lambda )  g,  h
\rangle = 0 \quad \quad \forall \lambda \in \Lambda $$
from which $c_\lambda = 0$ for all $\lambda $. 
\end{proof}

\section{New Characterizations  of Gabor Frames}

 We are now ready to  discuss  the new   criteria for Gabor frames. Precisely, for each of the operators associated to a Gabor
system, we will state a property that is equivalent to the frame
property. Conceptually, several of them are easier because they do not
involve the invertibility  or an inequality.  
The following  theorem provides a characterization of a Gabor frame $\gab $
with a test function in $g\in M^1(\rd )$ in terms of each of the associated
operators $C, D$ and their combinations.

\begin{tm}\label{char}
  Assume that $g\in M^1(\rd ), g\neq 0$. Then the following are
  equivalent:
  \begin{itemize}
  \item[(i)] $\gab $ is a frame for $\lrd $.  
\item[(ii)]  $S_{g,\latt } $ is invertible on
  $M^1(\rd )$.  
\item[(iii)] $S_{g,\latt } $ is invertible on
  $M^\infty (\rd )$.
\item[(iv)] $S_{g,\latt } $ is one-to-one  on
  $M^\infty (\rd )$. 
\item[(v)]  $C_{g,\latt } $   is one-to-one  from 
  $M^\infty (\rd )$ to $\ell ^\infty (\latt )$. 
\item[(vi)]  $D_{g,\latt } $  defined on  
  $\lone (\latt )$ has dense range in  $M^1(\rd )$.
 \item[(vii)] $D_{g,\latt } $   is surjective   from 
  $\lone (\latt )$ onto $M^1(\rd )$.
\item[(viii)]  $D_{g,\duall } $   is one-to-one    from 
  $\ell ^\infty  (\duall )$to $M^\infty (\rd )$.
 \item[(ix)]  $C_{g,\duall } $  defined on  
  $M^1(\rd )$ has dense range in  $\ell ^1(\latt )$.
\item[(x)]   $C_{g,\duall } $ is surjective from
  $M^1(\rd )$ onto   $\ell ^1(\latt )$. 
\item[(xi)] $G_{g,\duall }$ is invertible on $\lone (\duall )$.  
\item[(xii)]  $G_{g,\duall }$ is invertible on $\ell ^\infty  (\duall )$.  
\item[(xiii)]  $G_{g,\duall }$  is one-to-one  on $\ell
  ^\infty  (\duall )$.   
\item[(xiv)]  $\gabdual $ is a Riesz sequence in
  $\lrd $. 
 \end{itemize}
\end{tm}

The merit of Theorem~\ref{char} is its beauty and completeness. 
Conceptually it is simpler to verify the
injectivity of an operator than to prove its invertibility or
its surjectivity. To emphasize this point,  let us single out two conditions and
state separately the following simplified version of
Theorem~\ref{char}. The following was the original conjecture from
which Theorem~\ref{char} evolved. 

\begin{cor}\label{simple}
Assume that $g \in M^1(\rd )$ and $\latt $ is a lattice in $\rdd $
with adjoint lattice $\duall $. Then the following are equivalent:

\begin{itemize}
\item[(i)] The set $\gab $ is a frame for $\lrd $. 
\item[(ii)]The analysis operator $C_{g,\latt }$ is one-to-one from
  $M^\infty (\rd )$ to $\ell ^\infty (\Lambda )$. 
\item[(iii)] The synthesis operator $D_{g, \duall }$ is one-to-one from
  $\ell ^\infty (\Lambda )$ to $M^\infty (\rd )$.  
\end{itemize}
  \end{cor}

The verification of conditions (ii) or (iii) may be easier in some
cases. Each of these yields  the frame property of  $\gab $,
without involving  any  inequality or the inversion of an
operator. This insight seems to be completely new and is 
rather surprising.  In Section~5 we will  show how
Corollary~\ref{simple} can be used for counter-examples.

\vspace{3 mm}


\textbf{Diskussion of  Theorem~\ref{char}.}

 1.  The set of  conditions (viii) --- (xiv) is dual to the set of
conditions (i) --- (vii). The dual conditions are   obtained by
replacing the lattice  $\latt $ by its  adjoint 
lattice $\duall $ and interchanging the indices $1$ and $\infty $ or
interchanging the role of $C$ and $D$.  

The equivalence of (i) and (xiv) is the so-called \textrm{Ron-Shen
  duality} for Gabor frames. It is implicit in the work of
Rieffel~\cite{rieffel88}, was first obtained by
Janssen~\cite{janssen95} and then independently by Ron and
Shen~\cite{ron-shen97} and Daubechies, Landau, and
Landau~\cite{DLL95}. The version for arbitrary lattices is due to
Feichtinger and Kozek~\cite{feichtinger-kozek98}.  This duality is a fundamental 
principle in the \tfa\ of Gabor frames and is connected with deep
results in operator theory~\cite{rieffel88}. 
Theorem~\ref{char} extends  the duality theory of Gabor frames. 

2. It is well known that Gabor frames with windows in $M^1$ extend to
so-called Banach frames for the \modsp s. One way to formulate this
fact is the following: \emph{For $g\in M^1$ the following properties
  are equivalent. (a) The frame operator $S_{g,\latt}$ is invertible
  on $L^2$. (b) $S_{g,\latt }$ is invertible on $M^p$ for  some $p,
  1\leq p \leq \infty $. (c) $S_{g,\latt }$ is simultaneously
  invertible on all  $M^p$,    $   1\leq p \leq \infty $.   }
See ~\cite[Thm.~5.2]{fg97jfa} and \cite[Prop.~12.2.7]{book} and the theory of Gelfand
triples~\cite{feichtinger-kozek98}.  Whereas this theorem suggests that  all 
\modsp s play the same role,  Theorem~\ref{char} shows
that the situation in not entirely symmetric and that $M^\infty $
plays a distinguished role. On $M^\infty $ the frame operator $S_{g,
  \latt } $ is invertible, \fif\ it is one-to-one. This is far from
true on $L^2= M^2$. It is well-known that $S_{g,\latt } $ may be
one-to-one on $L^2$ without being invertible. See the examples in
Section~5.


\section{Proof of Theorem~\ref{char}}

Before we begin with the proof, we collect some facts from abstract
harmonic analysis. In particular, we review Wiener's Lemma for twisted
convolution and its role in the analysis of Gabor frames. 

  Some of the implications of  Theorem~\ref{char}  are quite easy to
prove,  but others require the full arsenal of \tfa . We will apply  
some of the deepest results  in \tfa . We only give a short
explanation of the results needed and place them in the  context in \tfa . Full
details are found in the cited  literature. 

\subsection{Twisted Convolution}

Recall first  \tfs s do not commute. If $\lambda =
(\lambda _1, \lambda _2) \in \rd \times \rd \simeq \rdd $ and $\mu  =
(\mu _1, \mu _2) \in \rdd $, then their composition is 
\begin{equation}
  \label{eq:33}
\pi (\lambda ) \pi (\mu ) = e^{-2\pi i \lambda _1  \cdot \mu _2 } \pi (\lambda + \mu )  \, .  
\end{equation}
The occurring phase factor defines  a  quadratic form $\sigma $ on
$\rdd$ by 
$$
\sigma (\lambda , \mu ) = \lambda _1  \cdot \mu _2 \, .
$$
Let $\bba , \bbb$ be two finite  sequences indexed by the lattice $\latt
\subseteq \rdd $. The twisted convolution $\bba \, \natural \, \bbb $
 is defined to be 
\begin{equation}
  \label{eq:31}
  (\bba \, \natural \, \bbb ) (\lambda ) = \sum _{\mu \in \latt }
  a_\lambda b_{\mu - \lambda } e^{2\pi i \sigma ( \lambda , \mu -
    \lambda )}
\end{equation}
Strictly speaking, the twisted convolution depends on the lattice
$\latt $ and we would have to write $\, \natural \, _\latt $. However,
since we use only a fixed lattice $\latt $ and its adjoint lattice
$\duall $, no confusion can arise, and  we
will omit the subscript. 

By Young's inequality the twisted convolution extends to certain $\ell
^p$-spaces; in particular we have $\|\bba \, \natural \, \bbb \|_p
\leq \|\bba \|_p \, \|\bbb \|_1$ for $\bba \in \ell ^p (\latt)$ and
$\bbb \in \ell ^1(\latt )$. 

In our context the fundamental property of twisted convolution is
\emph{Wiener's Lemma for twisted convolution}.

\begin{tm} \label{wltwist}
  Assume that $\bba \in \ell ^1(\latt )$ and that the (twisted
  convolution) operator $C_{\bba } \bbc = \bbc \, \natural \, \bba $
  is invertible on $\ell ^2(\latt )$. Then the inverse is  $C_{\bba
  }\inv = C_{\bbb }$ for the  unique  $\bbb \in \ell ^1 (\latt
  )$.  Consequently, $C_{\bba } $ is invertible simultaneously on
  \emph{all} $\ell ^p (\latt ), 1\leq p \leq \infty $. 
\end{tm}

For the proof and more general statements see~\cite{GL04}, for alternative proofs
see~\cite{GL04a,BCHL06}. 

\subsection{Twisted Convolution in Time-Frequency Analysis}
\label{twisttfa}

Twisted convolution arises in several contexts in \tfa , sometimes
naturally, sometimes hidden. 

(A) Consider series of \tfs s of the form $\pi (\bbc ) = \sum
_{\lambda \in \latt } c_\lambda \pi (\lambda )$. If $\bbc \in \ell
^\infty (\latt )$, then by Lemma~\ref{synthesis} $\pi (\bbc )$ is
bounded from the space of 
test functions  $M^1$ to the space of distributions 
$M^\infty $, so $\pi (\bbc ) $ is always well defined.   If $\bbc \in
\ell ^1(\latt )$, then $\pi (\bbc )$ is an 
\emph{absolutely convergent series of \tfs s} and is easily seen to be
bounded on each $M^p, 1\leq   p \leq \infty $. 

Let $\bba , \bbb \in \ell ^1(\latt )$, then the commutation
rule~\eqref{eq:33} implies that
\begin{equation}
  \label{eq:33b}
  \pi (\bba ) \pi (\bbb ) = \pi (\bba \twi \bbb ) \, .
\end{equation}
Furthermore, $\pi : \ell ^1(\latt ) \to \cB (\lrd )$ is a faithful
representation of the (involutive) Banach algebra $\ell ^1 (\latt )$
with respect to $\twi $ into the bounded operators on $\lrd
$~(\cite{GL04} or \cite{rieffel88}).  With some work, Wiener's Lemma
for twisted convolution can be transferred to the operator algebra
$\pi (\ell ^1(\latt ))$ of absolutely convergent \tfs s, see~\cite{GL04}.
As a result we state Wiener's Lemma for the so-called rotation algebra
$\pi (\ell ^1 (\Lambda ))$.

\begin{tm} \label{B}
  Let $S= \sum_{\lambda \in \Lambda } a_\lambda \pi (\lambda )= \pi
  (\bba )$. If  $\bba \in \ell ^1(\latt )$ and $S$ is invertible on
  $\lrd $, then $S\inv $ is again an absolutely convergent series of
  \tfs s $S\inv = \sum _{\lambda \in \Lambda } b_\lambda \pi (\lambda )= \pi
  (\bbb )$ for a unique $\bbb \in \ell ^1(\latt )$ satisfying $\bba
  \twi \bbb = \bbb \twi \bba = \delta $. As a consequence, $S$ is
  invertible simultaneously on all \modsp s $M^p, 1\leq p \leq \infty
  $. 
\end{tm}

\rem\ Although we do not need it here, we would like to point out an
interesting and deep relation between Gabor frames and operator
algebras and non-commutative geometry discovered by
Luef~\cite{lue06}. In the 
 language of operator algebras, $\pi (\ell 
^1(\latt ))$ is a \emph{rotation algebra} or
\emph{non-commutative torus}, and Wiener's Lemma is usually referred
to as the spectral invariance property or the spectral
permanence~\cite{arveson91,connes}. If the coefficient algebra $\ell ^1(\latt
)$ is replaced by the Frechet algebra $\cS (\latt )$ of rapidly
decaying sequences, the corresponding version of Theorem~\ref{B} is a
celebrated theorem of Connes about the spectral invariance of
\emph{smooth} non-commutative tori~\cite{connes80}.

\vspace{3 mm}

(B) \emph{The Gramian Operator.} Recall that 
$$
G_{\lambda , \mu} = \langle \pi (\mu  ) g, \pi (\lambda )
g\rangle  = e^{2\pi i \sigma (\mu , \lambda -\mu)} \langle g, \pi (\lambda -\mu)g\rangle
\, .
$$
Writing $a_\lambda =  \langle g, \pi (\lambda )g\rangle$, the
action of the Gramian operator $G$ can be written as a twisted
convolution:
\begin{eqnarray}
  (G\bbc )_\lambda &=& \sum _{\mu \in \latt } G_{\lambda \mu }c_\mu
  \notag \\
&=& \sum _{\mu \in \latt } c_\mu  \langle g, \pi (\lambda
-\mu)g\rangle\, e^{2\pi i \sigma (\mu , \lambda  -\mu )} \label{eq:32} \\
&=& (\bbc \twi \bba )_\lambda \, .   \notag    
\end{eqnarray}
This identity makes is plausible why Theorem~\ref{wltwist} enters the proof
of our main result. 

\vspace{2mm}

(C) \emph{Janssen's representation of the Gabor frame operator.}
The frame operator $S_{g,\latt }$ commutes with all \tfs s $\pi
(\lambda ) , \lambda \in \Lambda ,$ and thus belongs to the commutants
of $\pi (\Lambda )$. By definition of the adjoint lattice, this
commutant is spanned by the \tfs s $\{\pi (\mu ), \mu \in \duall \} $. It
is therefore plausible that $S_{g,\latt }$ can be represented by a sum
of \tfs s over the adjoint lattice $\duall $ in some sense. 
The precise statement is Janssen's representation~\cite{janssen95} of
the frame 
operator.

\begin{tm}
  If $g\in M^1$, then $S_{g,\latt } = \sum _{\mu \in \duall  } a_\mu
  \pi (\mu ) =  \pi (\bba )$ for some $\bba \in
  \ell ^1(\duall )$. The
coefficients are given explicitly by $a_\mu = s(\Lambda )\inv \langle
g, \pi (\mu ) g\rangle $, where $s(\Lambda ) = |\det A|$  is the size of   $\Lambda
= A\zdd $.
\end{tm}

If, in addition, $\gab $ is a frame, then $S_{g, \latt }$ is
invertible and Theorem~\ref{B} implies that $S_{g, \latt }$ is
invertible simultaneously on all \modsp s $M^p, 1\leq p \leq \infty
$.

\subsection{Plan of the Proof}

We will prove the following chain of
implications. 
First we  cover all conditions involving the lattice $\latt $ and show
how they imply condition (viii) involving the adjoint lattice $\duall
$. We will prove the implications 
\begin{equation}
  \label{eq:9}
  (i) \Rightarrow (ii) \Rightarrow (iii) \Rightarrow (iv) \Rightarrow
  (v) \Rightarrow (vi) \Rightarrow (viii) 
\end{equation}
and 
\begin{equation}
  \label{eq:10}
  (ii) \Rightarrow (vii) \Rightarrow (viii) \, .
\end{equation}

On the side of the adjoint lattice $\duall $ we will prove the
following implications:
\begin{equation}
  \label{eq:11}
  (viii) \Rightarrow (ix) \Rightarrow (x) \Rightarrow (xiv)
  \Rightarrow (i)
\end{equation}
and
\begin{equation}
  \label{eq:12}
(i) \Rightarrow  (xiv) \Rightarrow (xi) \Rightarrow (xii) \Rightarrow (xiii)
  \Rightarrow (viii) \, .
\end{equation}

 Let us start!

\subsection{Proof}
\label{sec}

\vspace{5 mm}

\textbf{(i) $\, \Leftrightarrow \, $ (ii). } The simultaneous
invertibility of the frame operator $S_{g,\Lambda }$ on $L^2 (\rd ),
M^1(\rd )$, and $M^\infty (\rd )$ is the main theorem of ~\cite{GL04};
For the special  lattices $\latt = \alpha \zd \times \beta \zd $ with
rational $\alpha \beta $ this fact  was already proved
in~\cite[Thm.~5.2]{fg97jfa}. See Section~\ref{twisttfa} for a discussion
of the context. 

\textbf{(ii) $\, \Leftrightarrow \, $ (iii). }  The frame operator
$S_{g,\latt } $ is self-adjoint. Hence it is invertible on $M^1(\rd )$ \fif\  it
is invertible on the dual space $M ^\infty (\rd ) $.

\textbf{(iii) \imp (iv) } is  obvious.

\textbf{(iv) \imp (v). } If $S_{g,\latt }= D_{g,\latt }C_{g,\latt }$ is one-to-one
on $M^\infty (\rd )$, then clearly  $ C_{g,\latt }$ must  be one-to-one
from $M^\infty (\rd )$ to $\ell ^\infty (\latt )$. 

\textbf{(v) $\, \Leftrightarrow \, $ (vi). }  The adjoint operator of
$D_{g,\latt} : \ell ^1(\Lambda ) \to M^1(\rd )$ is exactly $C_{g,\latt
} : M^\infty (\rd ) \to \ell ^\infty (\latt )$. So $C_{g,\latt }$ is
one-to-one, \fif\ its adjoint $D_{g,\latt }$ has dense range. 

\textbf{(vi), (vii) \imp (viii). } By assumption  $D_{g,
  \Lambda } \big(\ell ^1 (\latt )\big) $ is a dense subspace of
$M^1(\rd )$ (or equals $M^1(\rd )$). Then by Lemma~\ref{bound}(ii) the
finite linear combinations of the form $f = \sum _{\lambda  }
a_\lambda \pi (\lambda ) g\in M^1(\rd )$  span also a dense subspace
of $M^1$. 
  Now assume
that  
$$
\sum _{\mu \in \duall } c_\mu \pi (\mu ) g  = 0
$$
for some $\mathbf{c} \in \ell ^\infty (\duall )$ (as a distribution in
$M^\infty $).  Let  $f = \sum _{\lambda \in F } a_\lambda \pi (\lambda ) g\in
M^1(\rd )$ 
for some finite set $F\subseteq \latt $. 
Since $\pi (\lambda ), \lambda \in \latt $ and $\pi (\mu ), \mu \in 
\duall $ commute, we find  that
\begin{equation*}
   \sum _{\mu\in \duall } c_\mu \pi (\mu ) f =   \sum _{\mu\in \duall
   } c_\mu \pi (\mu ) \Big( \sum _{\lambda \in F} a_\lambda \pi (\lambda
  ) g\Big)  = \sum _{\lambda \in F}   a_\lambda \pi (\lambda
  ) \Big(\sum _{\mu\in \duall } c_\mu \pi (\mu ) g\Big) = 0 \, .   
\end{equation*}
 This calculation is justified  by Lemma~\ref{synthesis}. Thus  $ \sum
_\mu c_\mu \pi (\mu )   :M^1 \to M^\infty $ 
vanishes on a dense subspace of $M^1$, consequently, $ \sum
_\mu c_\mu \pi (\mu )=0$.  Since \tfs s are linearly independent by Proposition~\ref{linind}, it
follows that $\mathbf{c} = 0$. 

\textbf{(ii) \imp (vii).} If $S_{g,\latt }= D_{g,\latt }C_{g,\latt }$ is  a
bijection on   $M^1$, then $ D_{g,\latt }$ must be surjective
from $\ell ^1 (\latt )$ onto $M^1(\rd )$. 
 
\vspace{5mm}

\textbf{(viii) $\, \Leftrightarrow \, $  (ix). } Again, $C_{g,\duall
}: M^1 \to \ell ^1(\duall )$ 
has dense range, \fif\ its adjoint operator $D_{g,\duall }: \ell
^\infty (\duall ) \to M^\infty $ is one-to-one. 

 \textbf{(ix) \imp (x). } Assume that the analysis operator
 $C_{g,\duall } $ from  $M^1$ has dense range in $\ell ^1 (\duall
 )$. Then for any fixed  $\epsilon, 0\leq \epsilon <1$, there is a
 function $\varphi \in M^1(\rd )$ 
 such that $\|C_{g,\duall } \varphi  - \delta \|_{1} < \epsilon $,
 explicitly, 
 \begin{equation}
   \label{eq:16}
   \sum _{\mu \in \duall } |\langle \varphi , \pi (\mu ) g \rangle -
   \delta _{\mu ,0}| = |\langle \varphi,g\rangle -1| + \sum _{\mu \in
     \duall , \mu \neq 0} |\langle \varphi , \pi (\mu ) g | <\epsilon
   \, .
 \end{equation}
Now let $\Phi $ be the matrix defined by the  entries
\begin{equation}
  \label{eq:17}
  \Phi _{\mu \nu } = \langle \pi (\nu )\vf  , \pi (\mu )g \rangle \qquad
  \mu , \nu \in \duall \, .
\end{equation}
($\Phi $ is the ``cross Gramian'' of the Gabor
systems $\gabdual $ and $\cG (\vf, \duall )$.)
We will show that $\Phi $ is invertible by  applying  Schur's test to
$\Phi -\mathrm{I}$. First we estimate the operator norm on $\ell ^1$: 
\begin{eqnarray*}
\|\Phi - \mathrm{I}\|_{\ell ^1 \to \ell ^1} &=&   \sup _{\nu \in \duall } \sum _{\mu \in \duall } |\Phi _{\mu \nu } -
  \delta _{\mu \nu }| \\
& = &  \sup _{\nu \in \duall } \sum _{\mu \in
    \duall } | \langle \pi (\nu )\vf  , \pi (\mu )g \rangle - \delta
  _{\mu \nu }| \\
&=& \sup _{ \nu \in \duall } \Big( |\langle \varphi,g\rangle -1| +
\sum _{\mu \in      \duall , \mu \neq \nu} |\langle  \varphi
, \pi (\mu - \nu ) g\rangle  | \Big) \\ 
&=& |\langle \varphi,g\rangle -1| + \sum _{\mu \in
     \duall , \mu \neq 0} |\langle \varphi , \pi (\mu ) g \rangle | \\
&=& \|C_{g, \duall } \vf - \delta \|_1 <\epsilon
   \, .
\end{eqnarray*}

Since   $\|\Phi - \mathrm{I} \|_{\ell ^1 \to \ell ^1} < \epsilon < 1
$,  $\Phi $ is invertible on $\ell ^1(\duall )$. 

Now let $\mathbf{a} \in \ell ^1 (\duall )$ be arbitrary, then there
exists  $\mathbf{c} \in \ell ^1 (\duall )$, such that $\Phi \mathbf{c} =
\mathbf{a}$. Set $f = \sum _{\nu \in \duall } c_\nu \pi (\nu )\vf $,
then 
\begin{eqnarray}
  (C_{g, \duall } f)(\mu ) &=& \langle f, \pi (\mu ) g \rangle \notag \\
&=& \sum _{ \nu \in \duall } c_\nu \langle \pi (\nu )\vf , \pi (\mu )
g\rangle   \label{eq:19}
 \\
&=& (\Phi \mathbf{c})(\mu ) = a_\mu \, . \notag 
\end{eqnarray}
Thus $C_{g,\duall }$ is surjective from $M^1$ onto $\ell
^1(\duall )$. 

\textbf{(x) \imp (xiv). } If $C_{g,\duall }$ is surjective from
$M^1$ onto $\ell ^1(\duall )$, then there exists a function
$\gamma \in M^1(\rd )$, a so-called \emph{dual window}, such that 
\begin{equation}
  \label{eq:20}
  \langle \gamma , \pi (\mu ) g \rangle = \delta _{\mu ,0} \qquad
  \text{ for all } \, \mu \in \duall \, .
\end{equation}
(This is the so-called Wexler-Raz biorthogonality condition and known
to be equivalent to $\gab $ being a frame~\cite{DLL95,janssen95}. So
condition (x) implies (i) directly.) 
Consequently, we also have 
\begin{equation}
  \label{eq:21}
  \langle \pi (\nu )\gamma , \pi (\mu ) g \rangle = \delta _{\mu ,\nu} \qquad
  \text{ for all } \, \mu , \nu  \in \duall \, .
\end{equation}
If $f=\sum _{\mu \in \duall } c_\mu \pi (\mu )g= D_{g, \duall }
\mathbf{c}$, then, according to  \eqref{eq:21},  the coefficients are determined by 
\begin{equation}
  \label{eq:22}
  \langle f, \pi (\nu ) \gamma \rangle = \sum _{\mu \in \duall
  } c_\mu \langle \pi (\mu )g,  \pi (\nu ) \gamma \rangle = \sum _{\mu
    \in \duall   } c_\mu \delta _{\mu \nu } = c_\nu \, .
\end{equation}
In different notation,
\begin{equation}
  \label{eq:23}
  \mathbf{c} = C_{\gamma , \duall } f \, ,
\end{equation}
and by the boundedness of $C_{\gamma , \duall }$ 
(Lemma~\ref{bound}(i)) we find that
$$
\|\mathbf{c} \|_2 = \| C_{\gamma , \duall } f \|_2  \leq C \|\gamma
\|_{M^1}\,  \|f \|_2  =
C \|\gamma \|_{M_1} \, \|D_{\gamma , \duall } \mathbf{c} \|_2 \leq C^2
\|\gamma \|_{M^1}^2 \, \|\mathbf{c} \|_2 \, 
$$
for all $\mathbf{c} \in \ell ^2 (\duall )$. This inequality says that  
$\gabdual $ is a Riesz sequence, and we have proved (xiv).

\textbf{(xiv)  $\, \Leftrightarrow \, $  (i). } This equivalence is
known as  the duality principle for  Gabor frames and well known. See
the discussion in Section~3 and~\cite{janssen95,ron-shen97}.

\textbf{(xiv) \imp (xi), (xii). } We have already seen in~\eqref{eq:32}
that the action of the Gramian operator can be recast as a twisted
convolution operator. Set $a_\mu = \langle g, \pi (\mu ) g\rangle $
for $\mu \in \duall $, then $\bba \in \ell ^1(\duall )$ by~\eqref{eq:35} and
~\eqref{eq:32} expresses the action of $G$ as a twisted convolution: 
$$
G\bbc = \bbc \twi \bba \, .
$$

Now assume that $\gabdual $ is a Riesz sequence, then the Gramian
operator $G $ (and thus the twisted convolution operator)  is
invertible on $\ell ^2(\duall )$. By Wiener's Lemma
(Theorem~\ref{wltwist}) 
$G\inv $ is again a twisted convolution operator $G\inv \bbc = \bbc
\twi \bbb $ for some $\bbb \in \ell ^1(\duall )$. Consequently $G\inv
$ is invertible on all $\ell ^p (\duall )$, $1\leq p \leq \infty $, as
was to be proved.

\textbf{(xi) $\, \Leftrightarrow \, $  (xii). }  The matrix $G_{g,\duall
} $ is self-adjoint. Hence it is invertible on $\ell^1(\duall )$ \fif\ if it
is invertible on the dual space $\ell ^\infty (\duall ) $.

\textbf{(xii) \imp (xiii) } is obvious.

\textbf{(xiii) \imp (viii).} If $G_{g,\duall } = C_{g,\duall }
D_{g,\duall }$ is one-to-one on $\ell ^\infty (\duall )$, then clearly
$D_{g,\duall }$ must be one-to-one from $\ell ^\infty (\duall )$ to
$M^\infty (\rd )$. \\

We have now shown that all fourteen conditions are equivalent and the
proof is complete. \hfill  \qed

\section{ Examples and Further Topics}

We conclude with  several examples and some variations.

 In general, Theorem~\ref{char} will not make it
easier to verify that a 
particular Gabor system $\gab $ is a frame. However,
Theorem~\ref{char} can be used to show that  $\gab $
fails to be a frame  by constructing explicit
sequences in the kernel of $D_{g,\duall }$. We will do this in two
(known) cases. The benefit of Theorem~\ref{char} is again that no
inequalities have to be proved or disproved.

  (a) Let $g(t) = e^{-\pi t^2}$ be the Gaussian and $\Lambda = \duall
  = \bZ ^2$. It  was shown in ~\cite{bargmann71} that  the
  synthesis operator $D_{g, \bZ ^2}$ is one-to-one on $\ell ^2(\bZ
  ^2)$ and hence the linear combinations of \tfs s of the Gaussian are
  dense in $L^2( \bR )$. This fact was already  claimed by J.~von
  Neumann~\cite{von-neumann}, see also \cite{baastians80,janssen82}.  
By contrast, $D_{g, \bZ ^2}$ is \textbf{not} one-to-one on $\ell
^\infty (\bZ ^2)$. By Theorem~\ref{char} the Gabor system $\cG (g, \bZ
^2)$ cannot be a frame for $L^2 (\bR )$. An example in the kernel of
$D_{g, \bZ ^2} $ is the sequence $c_{kl}=(-1)^{k+l}$. One verifies
that 
$$
\sum _{k,l\in \bZ } (-1)^{k+l} e^{2\pi i k t } e^{-\pi (t-l)^2} = 0
$$
as a distribution in $M^\infty (\bR )$. In fact, $\mathrm{ker}\,
D_{g,\bZ ^2} = \bC \mathbf{c} $, the kernel  has dimension
$1$. Consequently, the synthesis operator $D_{g,\bZ ^2}$ is also  one-to-one
on $\ell ^p$ for $p< \infty $. This example shows that one may not
hope for a better result in Theorem~\ref{char}.

(b) If $g \in M^1(\bR  )$ satisfies the partition-of-unity condition 
$$
\sum _{k\in \zd } g(t-\gamma  k) = a \neq 0 \qquad \text{for all } \, t \in
\rd \, ,
$$
then the Gabor system $\cG (g, \Lambda )$ cannot be a frame for any
lattice of the form $\Lambda = \alpha \zd \times \frac{N}{\gamma }
\zd $, where $\alpha >0$ is arbitrary and 
$N  = 2, 3 , 4 , \dots $. This was proved in~\cite{delprete,GJKP04}. 
Theorem~\ref{char} offers a simple proof of this fact: we have to show
that $D_{g,\duall }$ possesses a non-trivial kernel in $\ell ^\infty
$. We first write
$$
D_{g,\duall } \mathbf{c} = \sum _{k,l \in \bZ } c_{kl} M_{\ell /\alpha
} T_{k/\beta } g = \sum _{k\in \bZ } m_k T_{k/\beta }g
$$ 
for a sequence of $\alpha $-periodic functions/distributions $m_k$
with bounded Fourier coefficients. Given $N\geq 2$,  choose a sequence
$m_k$ to be $N$-periodic, namely $m_{k+N} = m_k$ for all $k\in \bZ
$. Then
\begin{eqnarray*}
  \sum _{k\in \bZ } m_k T_{k/\beta }g &=& \sum _{j=0}^{N-1} \sum _{l
    \in \bZ } m_{j+lN}  T_{(j+lN)\alpha /N }g \\
&=& \sum _{j=0}^{N-1} m_j T_{j\alpha /N} \Big( \sum _{l
    \in \bZ }   T_{l\alpha }g \Big) = a \sum _{j=0}^{N-1} m_j \, .
\end{eqnarray*}
The latter sum vanishes with an appropriate choice of $m_j$, e.g.,
choose $m_{lN} = 1$ and $m_{1+lN} = -1$ and $m_{j+lN} = 0$ for $j=2,
\dots , N-1$ and $ l \in \bZ $. The corresponding sequence $\mathbf{c}$ is $c_{j+mN, l}
= (-1)^j \delta _{l,0}$ for $j=0,1$ and $l,m \in \bZ $ and $c_{kl}= 0$
otherwise. Of course, a similar statement can be formulated in
dimension $d>1$. 

If $g$ is the  finite or infinite convolution product $$g= \,^* \prod _{j=1}^\infty \chi
_{[-a_j/2, a_j/2]} = \chi
_{[-a_1/2, a_1/2]} \ast  \chi
_{[-a_2/2, a_2/2]}\ast \chi
_{[-a_3/2, a_3/2]}\ast \dots $$ for $a_j > 0$ and $\sum a_j < \infty $, then
$\sum _{k\in \bZ } g(x-a_j k) = b_j \neq 0$ for each
$j$. Consequently, for any lattice of the form $\Lambda = \alpha \bZ
\times \beta \bZ $ with $\alpha >0 $ arbitrary and $\beta = N/a_j$,
$j\in \bN $ and $N=2,3, \dots $, the Gabor system 
$\cG (g, \latt )$ cannot be a frame for $L^2$. This example shows that
the pattern of excluded frequency parameters $\beta $ can be made
arbitrarily complicated. 

\vspace{3 mm}

 \textbf{A direct proof of the implication \textbf{(xiii)
  \imp (xi), (xii)}}. We note this implication can be proved directly
by using  an observation of Y.~Choi in a more
general context. The sequence
space $\ell ^1 (\duall )$ is a Banach algebra under twisted
convolution $\, \natural \, $, and its dual is $\ell ^\infty (\duall
)$. Choi~\cite[Lemma~2]{choi06} observed that if the (twisted) convolution
operator $ \mathbf{c} \to \mathbf{c} \, \natural \, \mathbf{s} $ is
one-to-one on $\ell ^\infty $, then the adjoint operator, which is
again a convolution with a sequence $\mathbf{s}^*$ is surjective on
$\ell ^1 (\duall )$. In the case of the Gramian, the twisted
convolution is self-adjoint $\mathbf{s} = \mathbf{s}^*$, and so the
injectivity of twisted convolution by $\mathbf{s}$ on $\ell ^\infty $
implies its invertibility on $\ell ^1$.  

\vspace{5 mm}

\textbf{ Structure of $\mathrm{ker} \, D_{\duall }\ell ^\infty $.} If
$D_{g, 
   \duall } $ is non-trivial, then the kernel possesses special
 invariance properties that might be useful in order to disprove that
 a Gabor system is a frame. 

 \begin{lemma}
   The subspace $\mathrm{ker}\, D_{g,\duall }$ is w$^*$-closed and  an $\ell ^1$-module
   under twisted convolution. Thus, if $\bba \in \ell ^1(\duall )$ and
   $\bbc \in \mathrm{ker}\, D_{g, \duall }$, then $\bba \twi \bbc \in
   \mathrm{ker}\, D_{g, \duall }$. 
 \end{lemma}

 \begin{proof}
   Let $\bba \in \ell ^1(\duall )$ and $D_{g,\duall } \bbc = \sum
   _{\mu \in \duall } c_\mu \pi (\mu ) g = 0$ in $M^\infty $. Then as
   in~\eqref{eq:33b} we find that
   \begin{eqnarray*}
0&=&      \pi (\bba ) \sum _{\mu \in \duall } c_\mu \pi (\mu )g \\
&=&  \pi      (\bba ) \pi (\bbc ) g = \pi (\bba \twi \bbc ) g \\
&=& D_{g, \duall } (\bba \twi \bbc ) \, ,
   \end{eqnarray*}
and thus $\bba \twi \bbc \in \mathrm{ker}\, D_{g, \duall }$. The
w$^*$-closedness is clear.  
 \end{proof}

The module property of $\mathrm{ker}\, D_{g,\latt} $ suggests  the
definition of an index for a Gabor system $\gab $. 

\begin{definition}
  The \emph{index}  of the Gabor system
  $\gab $, denoted by $\mathrm{ind}\, (g,\latt )$,  is the smallest cardinality of a set of module generators of
  $\mathrm{ker}\, D_{g,\latt }$. We write $\mathrm{ind}\, (g,\latt ) =
  N$, if there exist $N$ sequences  $e_j \in \mathrm{ker} \, D_{g, \latt
  }$, such that every $\bbc \in \mathrm{ker} \, D_{g, \latt
  }$ can be written as 
$$
\bbc = \sum _{j=1}^N \bba _j \twi e_j$$
for some $\bba _j \in \ell ^1(\latt )$ and there is no set of smaller
cardinality with this property. 
\end{definition}

With this definition we can recast the equivalence \textbf{(i) $\,
  \Leftrightarrow \, $ (viii)} as follows.

\begin{tm} \label{index}
  The set $\gab $ is a frame for $\lrd $ \fif\ $\mathrm{ind} \,
  (g,\duall ) = 0$. 
\end{tm}

At this point it is not clear what is the real significance of the
index and what it says about  rotation algebras and related
objects. Here is what we know so far:

(a) If the \tfs s $\{\pi (\lambda ), \lambda \in \latt \}$ commute with
each other, e.g. for the lattice $\Lambda = \zdd $, then $\mathrm{ker} \, D_{g, \latt
  }$ is a translation-invariant weak$^*$-closed subspace of $\ell
  ^\infty $. By a Tauberian theorem of Wiener, $\mathrm{ker} \, D_{g, \latt
  }$  contains a sequence of the form $e_\xi (k) = e^{2\pi i \xi \cdot
    k}, k\in \zdd $. Furthermore, $\mathrm{ind}\, (g, \zdd ) =
  \mathrm{card}\, \{\xi \in \bT ^{2d}: e_\xi \in \mathrm{ker} \, D_{g, \latt
  }\}$. 

(b) Since every lattice is of the form $\Lambda = A\zdd $ for $A\in
\mathrm{GL}(2d, \bR )$, we may endow the set of lattices with the
topology of $\mathrm{GL}(2d, \bR )$. By a theorem of Feichtinger and
Kaiblinger~\cite{FK02}, the set of $g\in M^1$ and lattice $\latt $ for which $\gab
$ is a frame is open in $M^1 \times \mathrm{GL}(2d, \bR )$. Thus by
Theorem~\ref{index} the preimage $\mathrm{ind}\inv (\{0\})$ is open in
$M^1 \times \mathrm{GL}(2d, \bR )$. This suggests that the index is a
continuous function on $M^1 \times \mathrm{GL}(2d, \bR )$.

\vspace{3 mm}

\textbf{Other Types of Gabor Sets.} The  characterization of
Theorem~\ref{char} can be extended to Gabor sets with several basis
functions, so-called multi-window Gabor sets,  and to vector-valued
Gabor frames, so-called superframes. The set $\bigcup _{j=1}^n \cG
(g_j, \Lambda )$   is a frame, if for some constants $A,B>0$,
$$
A\|f\|_2^2 \leq \sum _{\lambda \in \Lambda } \sum _{j=1}^n |\langle f,
\pi (\lambda ) g_j\rangle | ^2 \leq B \|f\|_2^2 \qquad \forall f\in \lrd
\, .
$$ 
The associated frame operator is $S= \sum _{j=1}^n S_{g_j, \Lambda
}$. The same arguments as in the scalar case of Theorem~\ref{char}
show that, for $g_j \in M^1$, $S$ is invertible on $\lrd $ \fif\ $S$
is invertible on 
$M^p$ for some/all $p, 1\leq p \leq \infty $. The duality can be
expressed as follows: let $\mathbf{g} = (g_1, g_2, \dots , g_n) $ be
the vector with component functions $g_j $ and $\pi (\lambda )
\mathbf{g} = (\pi (\lambda )g_1, \dots , \pi (\lambda )g_n)$ (action
of $\pi $ is componentwise).  Then $\bigcup _{j=1}^n \cG
(g_j, \Lambda )$ is a frame for $\lrd $, \fif\ $\cG (\mathbf{g}, \duall ):= \{ \pi
(\mu ) \mathbf{g}: \mu \in 
\duall \}$ is a Riesz sequence in the space $L^2(\rd , \bC ^n)$ of
vector-valued functions. The Gramian of $\cG (\mathbf{g}, \duall )$
has the entries $G_{\mu ,\mu '} = \langle \pi (\mu ' ) \mathbf{g}, \pi
(\mu ) \mathbf{g}\rangle = \sum _{j=1}^n  \langle \pi (\mu ' ) g_j, \pi
(\mu ) g_j\rangle$, hence the same arguments as in Section~\ref{sec}
show that $G$ is invertible on $\ell ^2(\duall )$ \fif\ it is
invertible on $\ell ^p (\duall ), 1\leq p \leq \infty $. The
conditions involving the coefficient and synthesis operator need a bit
more care, because in this case they are different for $\Lambda $ and
$\duall $. Specifically, $D_{\mathbf{g},\Lambda }$ maps $\ell ^p
(\Lambda , \bC ^n)$ to $M^p$ by $D_{\mathbf{g}, \Lambda }(\bbc _j) =
\sum _{\lambda \in \Lambda } \sum _{j=1}^n c_{\lambda , j} \pi
(\lambda )g_j$, whereas $D_{\mathbf{g}, \duall }$ maps $\ell ^p
(\duall )$ to $M^p(\rd, \bC ^n)$ by  $D_{\mathbf{g}, \duall }(\bbc ) =
\sum _{\mu  \in \duall } c_{\mu} \pi
(\mu )\mathbf{g}$. With these precautions all conditions of
Theorem~\ref{char} can be formulated adequately and yield a
characterization of ``multi-window Gabor frames'' without
inequalities. We will return to this topic. 

Finally, we mention that Theorem~\ref{char} carries over without
change  to characterize Gabor frames in arbitrary locally compact
abelian groups. 

\def\cprime{$'$} \def\cprime{$'$}


\begin{thebibliography}{10}

\bibitem{arveson91}
W.~Arveson.
\newblock Discretized {CCR} algebras.
\newblock {\em J. Operator Theory}, 26(2):225--239, 1991.

\bibitem{baastians80}
M.~J. Baastians.
\newblock Gabor's expansion of a signal into {G}aussian elementary signals.
\newblock {\em Proc. IEEE}, 68:538--539, April 1980.

\bibitem{BCHL06}
R.~Balan, P.~G. Casazza, C.~Heil, and Z.~Landau.
\newblock Density, overcompleteness, and localization of frames. {II}. {G}abor
  systems.
\newblock {\em J. Fourier Anal. Appl.}, 12(3):309--344, 2006.

\bibitem{bargmann71}
V.~Bargmann, P.~Butera, L.~Girardello, and J.~R. Klauder.
\newblock On the completeness of coherent states.
\newblock {\em Rep. Math. Phys.}, 2:221--228, 1971.

\bibitem{choi06}
Y.~Choi.
\newblock Injective convolution operators on ${\ell}^{\infty}(\gamma)$ are
  surjective.
\newblock {\em Preprint}, 2006.

\bibitem{chr03}
O.~Christensen.
\newblock {\em An introduction to frames and {R}iesz bases}.
\newblock Applied and Numerical Harmonic Analysis. Birkh\"auser Boston Inc.,
  Boston, MA, 2003.

\bibitem{connes80}
A.~Connes.
\newblock {$C\sp{\ast} $} alg\`ebres et g\'eom\'etrie diff\'erentielle.
\newblock {\em C. R. Acad. Sci. Paris S\'er. A-B}, 290(13):A599--A604, 1980.

\bibitem{connes}
A.~Connes.
\newblock {\em Noncommutative geometry}.
\newblock Academic Press Inc., San Diego, CA, 1994.


\bibitem{daubechies92}
I.~Daubechies.
\newblock {\em Ten lectures on wavelets}.
\newblock Society for Industrial and Applied Mathematics (SIAM), Philadelphia,
  PA, 1992.

\bibitem{DLL95}
I.~Daubechies, H.~J. Landau, and Z.~Landau.
\newblock Gabor time-frequency lattices and the {W}exler-{R}az identity.
\newblock {\em J. Fourier Anal. Appl.}, 1(4):437--478, 1995.

\bibitem{delprete}
V.~Del~Prete.
\newblock Estimates, decay properties, and computation of the dual function for
  {G}abor frames.
\newblock {\em J. Fourier Anal. Appl.}, 5(6):545--562, 1999.

\bibitem{feichtinger80cras}
H.~G. Feichtinger.
\newblock Un espace de {B}anach de distributions temp\'er\'ees sur les groupes
  localement compacts ab\'eliens.
\newblock {\em C. R. Acad. Sci. Paris S\'er. A-B}, 290(17):A791--A794, 1980.

\bibitem{feichtinger81}
H.~G. Feichtinger.
\newblock On a new {S}egal algebra.
\newblock {\em Monatsh. Math.}, 92(4):269--289, 1981.

\bibitem{feichtinger89}
H.~G. Feichtinger.
\newblock Atomic characterizations of modulation spaces through {G}abor-type
  representations.
\newblock In {\em Proc. Conf. Constructive Function Theory, Edmonton, July
  1986}, pages 113--126, 1989.

\bibitem{feiSTSIP}
H.~G. Feichtinger.
\newblock Modulation spaces: looking back and ahead.
\newblock {\em Sampl. Theory Signal Image Process.}, 5(2):109--140, 2006.

\bibitem{fg97jfa}
H.~G. Feichtinger and K.~Gr{\"o}chenig.
\newblock Gabor frames and time-frequency analysis of distributions.
\newblock {\em J. Functional Anal.}, 146(2):464--495, 1997.

\bibitem{FK02}
H.~G. Feichtinger and N.~Kaiblinger.
\newblock Varying the time-frequency lattice of {G}abor frames.
\newblock {\em Trans. Amer. Math. Soc.}, 356(5):2001--2023 (electronic), 2004.

\bibitem{feichtinger-kozek98}
H.~G. Feichtinger and W.~Kozek.
\newblock Quantization of {T}{F} lattice-invariant operators on elementary
  {L}{C}{A} groups.
\newblock In {\em Gabor analysis and algorithms}, pages 233--266. Birkh\"auser
  Boston, Boston, MA, 1998.

\bibitem{FZ98}
H.~G. Feichtinger and G.~Zimmermann.
\newblock A {B}anach space of test functions for {G}abor analysis.
\newblock In {\em Gabor analysis and algorithms}, pages 123--170. Birkh\"auser
  Boston, Boston, MA, 1998.

\bibitem{book}
K.~Gr{\"o}chenig.
\newblock {\em Foundations of time-frequency analysis}.
\newblock Birkh\"auser Boston Inc., Boston, MA, 2001.

\bibitem{GJKP04}
K.~Gr{\"o}chenig, A.~J. E.~M. Janssen, N.~Kaiblinger, and G.~E. Pfander.
\newblock Note on {$B$}-splines, wavelet scaling functions, and {G}abor frames.
\newblock {\em IEEE Trans. Inform. Theory}, 49(12):3318--3320, 2003.

\bibitem{GL04}
K.~Gr{\"o}chenig and M.~Leinert.
\newblock Wiener's lemma for twisted convolution and {G}abor frames.
\newblock {\em J. Amer. Math. Soc.}, 17:1--18, 2004.

\bibitem{GL04a}
K.~Gr{\"o}chenig and M.~Leinert.
\newblock Symmetry and inverse-closedness of matrix algebras and functional
  calculus for infinite matrices.
\newblock {\em Trans. Amer. Math. Soc.}, 358(6):2695--2711 (electronic), 2006.

\bibitem{janssen82}
A.~J. E.~M. Janssen.
\newblock Bargmann transform, {Z}ak transform, and coherent states.
\newblock {\em J. Math. Phys.}, 23(5):720--731, 1982.

\bibitem{janssen95}
A.~J. E.~M. Janssen.
\newblock Duality and biorthogonality for {W}eyl-{H}eisenberg frames.
\newblock {\em J. Fourier Anal. Appl.}, 1(4):403--436, 1995.

\bibitem{lue06}
F.~{L}uef.
\newblock {O}n spectral invariance of non-commutative tori.
\newblock In {\em {O}perator {T}heory, {O}perator {A}lgebras, and
  {A}pplications}, volume 414, pages 131--146. {A}merican {M}athematical
  {S}ociety, 2006.

\bibitem{von-neumann}
J.~v. Neumann.
\newblock {\em Mathematische {G}rundlagen der {Q}uantenmechanik}.
\newblock Springer, Berlin, 1932.
\newblock English translation: ``Mathematical foundations of quantum
  mechanics,'' Princeton Univ. Press, 1955.

\bibitem{rieffel88}
M.~A. Rieffel.
\newblock Projective modules over higher-dimensional noncommutative tori.
\newblock {\em Canad. J. Math.}, 40(2):257--338, 1988.

\bibitem{ron-shen97}
A.~Ron and Z.~Shen.
\newblock Weyl--{H}eisenberg frames and {R}iesz bases in ${L}\sb 2({\bR } \sp
  d)$.
\newblock {\em Duke Math. J.}, 89(2):237--282, 1997.

\end{thebibliography}

\end{document}